\documentclass[12pt]{amsart}
\usepackage{amsmath,amsthm,amsfonts,amssymb,eucal}
\renewcommand {\a}{ \alpha }

\newcommand{\vare}{\varepsilon}

\newcommand{\varf}{\varphi}
\renewcommand{\d}{\delta}

\newcommand{\D}{\Delta}

\newcommand{\Sg}{\Sigma}
\renewcommand{\l}{\lambda}

\newcommand{\z}{\zeta}

\newcommand{\vart}{\vartheta}

\newcommand{\om}{\omega}

\newcommand{\R}{ \mathbb R}

\newcommand{\N}{ \mathbb N}

\newcommand{\Sq}{ \mathbb S}

\newcommand{\CB}{\mathcal B}

\newcommand{\CF}{\mathcal F}
\newcommand{\CG}{\mathcal G}
\newcommand{\CH}{\mathcal H}

\newcommand{\CP}{\mathcal P}

\newcommand {\bx}{\mathbf x}
\newcommand {\by}{\mathbf y}
\newcommand {\bb}{\mathbf b}

\newcommand {\BB}{\mathbf B}

\newcommand {\BF}{\mathbf F}

\newcommand {\BH}{\mathbf H}

\newcommand {\BM}{\mathbf M}

\newcommand {\BT}{\mathbf T}

\newcommand{\gz}{\mathfrak z}

\newcommand{\wt}{\widetilde}
\newcommand{\wh}{\widehat}

\DeclareMathOperator{\re}{Re}

\newtheorem{thm}{Theorem}[section]

\newtheorem{prop}[thm]{Proposition}

\theoremstyle{definition}

\theoremstyle{remark}

\numberwithin{equation}{section}

\newcommand{\thmref}[1]{Theorem~\ref{#1}}

\newcommand{\sh}{Schr\"odinger }
\newcommand{\vs}{\vskip0.2cm}

\newcommand{\loc}{\rm{loc}}
\newcommand{\w}{\infty}

\newcommand{\semi}{\rm{semi}}
\newcommand{\W}{\rm{Weyl}}
\newcommand{\rad}{\rm{rad}}
\newcommand{\nrad}{\rm{nrad}}

\begin{document}

\title[Spectral estimates for Schr\"odinger operators in $2D$]{On spectral estimates for two-dimensional Schr\"odinger operators}

\author[A. Laptev]{A. Laptev}
\address{Department of Mathematics\\ Imperial College London\\ Huxley Building\\
180 Queen's Gate\\ London SW7 2AZ, UK}
\email{a.laptev@imperial.ac.uk}
\author[M. Solomyak]{M. Solomyak}
\address{Department of Mathematics
\\ Weizmann Institute\\ Rehovot\\ Israel}
\email{michail.solomyak@weizmann.ac.il}

\subjclass[2010] {35J10; 35P20}
\keywords{\sh operator on $\R^2$; bound states; spectral estimates}

\begin{abstract}
For a two-dimensional \sh operator $\BH_{\a V}=-\D-\a V,\ V\ge 0,$
we study the behavior of the number $N_-(\BH_{\a V})$ of its negative eigenvalues (bound states), as the coupling parameter $\a$ tends to infinity. A wide class of potentials is described, for which $N_-(\BH_{\a V})$ has the semi-classical behavior, i.e., $N_-(\BH_{\a V})=O(\a)$. For the potentials from this class, the necessary and sufficient condition is found for the validity of the Weyl asymptotic law. \end{abstract}

\maketitle

\section{Introduction}\label{intro}
\subsection{Preliminaries}\label{prel}
Let $\BH_{\a V}$ be a \sh operator
\begin{equation}\label{sh}
    \BH_{\a V}=-\D-\a V
\end{equation}
on $\R^2$. We suppose that $V\ge 0$, and $\a>0$ is the coupling constant.
We write $N_-(\BH_{\a V})$ for the number of negative eigenvalues of $\BH_{\a V}$, counted with multiplicities:
\begin{equation*}
    N_-(\BH_{\a V})=\#\{j\in\N:\l_j(\BH_{\a V})<0\}.
\end{equation*}
As it is well known, the lowest possible (semi-classical) rate of growth of this function is
\begin{equation}\label{semic}
    N_-(\BH_{\a V})=O(\a),\qquad\a\to\infty.
\end{equation}
This agrees with the Weyl-type asymptotic formula
\begin{equation}\label{weyl}
    \lim_{\a\to\infty}\a^{-1}N_-(\BH_{\a V})=\frac1{4\pi}\int_{\R^2}Vdx
\end{equation}
that is satisfied if the potential behaves fine enough.

The exhaustive description of the classes of potentials on $\R^2$, such that \eqref{semic} or \eqref{weyl} is satisfied, is unknown till now. This is in contrast with the case of dimensions $d>2$, where the celebrated Cwikel-Lieb-Rozenblum estimate describes the class of potentials, for which both the estimate $N_-(\BH_{\a V})=O(\a^{d/2})$ and the Weyl asymptotic formula hold true.

In the forthcoming discussion, $\CP_{\semi}$ stands for the class of all potentials $V\ge0$ on $\R^2$, such that \eqref{semic} is satisfied, and $\CP_{\W}$ stands for the class of all such potentials that the asymptotics
\eqref{weyl} holds true.
It is clear that
\begin{equation}\label{semi-w}
   \CP_{\W}\subset\CP_{\semi}.
\end{equation}

The first results describing wide classes of potentials $V\in\CP_{\semi}$ were obtained in \cite{S-dim2} and \cite{BL}. In the latter paper, this was done also for the class $\CP_{\W}$. In particular, it was shown there that the inclusion in \eqref{semi-w}
is proper. What is more, in \cite{BL} the general nature of potentials
$V\in\CP_{\semi}\setminus\CP_{\W}$ was explained.

Some further estimates guaranteeing $V\in\CP_{\semi}$ were obtained in the recent paper \cite{GN}. We would like to mention also the paper \cite{mv} whose
authors have obtained some new results that give for $N_-(\BH_{\a V})$ the
order of growth larger than $O(\a)$.\vs

In the papers \cite{Lap, cha} the important case of radial potentials, $V(x)=F(|x|)$, was analyzed. For such potentials in \cite{cha} an integral estimate for $N_-(\BH_{\a V})$ was obtained  guaranteeing the inclusion $V\in \CP_{\semi}$ (actually, it guarantees also that $V\in \CP_{\W}$). This result was strengthened in the recent paper \cite{LS} where for the radial potentials the necessary and sufficient conditions for $V\in\CP_{\semi}$ and for $V\in \CP_{\W}$ were established.

In the present paper we return to the study of general (that is, not necessarily radial) potentials. We obtain an estimate that covers the main results of \cite{BL,GN}. It does not cover the estimate obtained in \cite{S-dim2}, however it has an important advantage compared with the latter: it does not use the intricate Orlicz norms appearing in \cite{S-dim2}.

In the papers \cite{Lap, cha} the important case of radial potentials, $V(x)=F(|x|)$, was analyzed. For such potentials in \cite{cha} an integral estimate for $N_-(\BH_{\a V})$ was obtained  guaranteeing the inclusion $V\in \CP_{\semi}$. Actually, it guarantees also that $V\in \CP_{\W}$. This result was strengthened in the recent paper \cite{LS} where for the radial potentials the necessary and sufficient conditions for $V\in\CP_{\semi}$ and for $V\in \CP_{\W}$ were established.

In the present paper we return to the study of general (that is, not necessarily radial) potentials. We obtain an estimate that covers the main results of \cite{BL}. It does not cover the estimate obtained in \cite{S-dim2},
however it has an important advantage compared with the latter: it does not use the intricate Orlicz norms appearing in \cite{S-dim2}.

\subsection{Formulation of the main result.}\label{basres} Below $(r,\vart)$ stand for the polar coordinates in $\R^2$, and $\Sq$ stands for the unit circle $r=1$. Given a function $V$, such that $V(r,\cdot)\in L_1(\Sq)$ for almost all $r>0$, we introduce its radial and non-radial parts
\begin{gather*}
    V_{\rad}(r)=\frac1{2\pi}\int_{\Sq} V(r,\vart)d\vart;\qquad 
    V_{\nrad}(r,\vart)=V(r,\vart)-V_{\rad}(r).
\end{gather*}

In our result the conditions will be imposed separately on
the radial and on the non-radial parts of a given potential $V$.
For handling the radial part, we need some auxiliary operator family on the real line, of the form
\begin{equation}\label{dim1}
    (\BM_{\a G}\varf)(t)=-\varf''(t)-\a G(t)\varf(t),\qquad \varf(0)=0,
\end{equation}
with the "effective potential"
\begin{equation}\label{effpot}
    G(t)=G_V(t)=e^{2|t|}V_{\rad}(e^t).
\end{equation}
Due to the condition $\varf(0)=0$ in \eqref{dim1}, for every $\a$ the operator
$\BM_{\a G}$ is the direct sum of two operators, each acting on the half-line. The sharp spectral estimates for $\BM_{\a G}$ can be given in terms of the number sequence (see Eq. (1.13) in \cite{LS})
\begin{equation}\label{sum}
   \wh{ \gz}(G)=\{\wh{\z_j}(G)\}_{j\ge0}:\qquad
     \wh{\z_0}(G)=\int_{D_0}G(t)dt, \quad
   \wh{\z_j}(G)=\int_{|t|\in D_j}|t|G(t)dt\quad  (j\in\N)
   \end{equation}
where $D_0=(-1,1)$ and $D_j=(e^{j-1},e^j)$ for $j\in\N$. For our purposes, it is convenient to express properties of this sequence in terms of the "weak $\ell_q$-spaces" $\ell_{q,\w}$.
Actually,  in the main body of this paper we deal only with $q=1$, and below we remind the definition of $\ell_{1,\w}$. The definition of the weak $\ell_q$-spaces with $q\neq1$ can be found, e.g., in \cite{BL}, Subsection 1.4.

 Given a sequence of real numbers $\bx=\{x_j\}_{j\in\N}$, such that $x_j\to 0$, we denote \begin{equation*}
    n_+(\vare,\bx)=\#\{j:|x_j|>\vare\},\qquad\vare>0.
 \end{equation*}
 The sequence $\bx$ belongs to $\ell_{1,\w}$, if
\begin{equation*}
    \|\bx\|_{1,\w}:=\sup_{\vare>0}(\vare\, n_+(\vare,\bx))<\infty.
\end{equation*}
This is a linear space, and the functional $\|\cdot\|_{1,\w}$ defines a quasinorm in it. The latter means that, instead of the standard triangle
inequality, this functional meets a weaker property:
\[  \|\bx+\by\|_{1,\w}\le c\bigl(\|\bx\|_{1,\w}+\|\by\|_{1,\w}\bigr),\]
with some constant $c>1$ that does not depend on the sequences $\bx,\by$. This quasinorm defines a topology in $\ell_{1,\w}$; there is no norm compatible with this topology.

 The space $\ell_{1,\w}$ is non-separable. Consider its closed subspace $\ell_{1,\w}^\circ$ in which the sequences $\bx$
with only a finitely many non-zero terms form a dense subset. This subspace is separable, and its elements are characterized by the property
\[\bx\in\ell_{1,\w}^\circ \ \Longleftrightarrow\ \vare\ n_+(\vare,\bx)\to0,\qquad \vare\to0.\]
The (non-linear) functionals
\begin{equation}\label{Delta}
    \D_1(\bx)=\limsup_{\vare\to0}(\vare\ n_+(\vare,\bx)),\qquad
    \d_1(\bx)=\liminf_{\vare\to0}(\vare\ n_+(\vare,\bx))
\end{equation}
are well-defined on the space $\ell_{1,\w}$, and
\[ \d_1(\bx)\le \D_1(\bx)\le \|\bx\|_{1,\w}.\]
It is clear that $\ell_{1,\w}^\circ=\{\bx\in\ell_{1,\w}: \D_1(\bx)=0\}$.\vs

The conditions on $V_{\nrad}$ will be given in terms of the space $L_1\left(\R_+,\, L_p(\Sq)\right)$, with an arbitrarily chosen $p>1$. This is the function space on $\R^2$, with the following norm:
\begin{equation}\label{L1Lp}
    \|f\|_{L_1\left(\R_+,\, L_p(\Sq)\right)}=
    \int_{\R_+}\left(\int_{\Sq}|f(r,\vart)|^pd\vart\right)^{1/p}rdr.
\end{equation}
This is a separable Banach space, and the bounded functions whose support is a compact subset in $\R^2\setminus\{0\}$ are dense in it. The space $L_1\left(\R_+,\, L_p(\Sq)\right)$ was used in the
paper \cite{LN}, and its results are one of the basic tools in our proof below.\vs

Here is the main result of the paper.

\begin{thm}\label{main}
Let a potential $V\ge0$ be such that $\wh{ \gz}(G_V)\in\ell_{1,\w}$, and
\begin{equation}\label{nrad-cond}
    V_{\nrad}\in L_1\left(\R_+,\, L_p(\Sq)\right)\qquad{\rm{with\ some}}\ p>1.
\end{equation}
Then $V\in\CP_{\semi}$, and the estimate is satisfied
 \begin{equation}\label{estim}
    N_-(\BH_{\a V})\le 1+C(p)\a\left(\| V_{\nrad}\|_{L_1\left(\R_+,\, L_p(\Sq)\right)}+
    \|\wh{ \gz}(G_V)\|_{\ell_{1,\w}}\right).
 \end{equation}

Moreover, the following equalities hold true:
\begin{equation}\label{as2}
   \begin{cases}
    &\limsup\limits_{\a\to\w}\a^{-1}N_-(\BH_{\a V})=\frac1{4\pi}\int Vdx+
    \limsup\limits_{\a\to\w}\a^{-1}N_-(\BM_{\a G_V}),\\
    &\liminf\limits_{\a\to\w}\a^{-1}N_-(\BH_{\a V})=\frac1{4\pi}\int Vdx+
    \liminf\limits_{\a\to\w}\a^{-1}N_-(\BM_{\a G_V}).
    \end{cases}
\end{equation}

 In particular, under the assumption \eqref{nrad-cond} the condition $\wh{ \gz}(G_V)\in\ell_{1,\w}^\circ$ is necessary and sufficient for
 $V\in\CP_{\W}$.
\end{thm}

In \eqref{as2}, and later on, the integral with no domain specified always means $\int_{\R^2}$.\vs

The formula \eqref{as2}, and especially, its proof
in Subsection \ref{prmain-as}, show that, in a certain sense, the parts $V_{\rad}$ and $V_{\nrad}$ contribute to the asymptotic behavior of $N_-(\BH_{\a V})$ independently. It may also happen that the contribution of $V_{\rad}$ is stronger than that of $V_{\nrad}$, and "screens" the latter.
This situation is described by the following statement, that complements our main theorem.
\begin{prop}\label{add}
Let a potential $V\ge 0$ be such that $\wh{ \gz}(G_V)\in\ell_{q,\w}$ with some $q>1$, and \eqref{nrad-cond} is satisfied. Then
\begin{equation*}\label{as2q}
   \begin{cases}
    &\limsup\limits_{\a\to\w}\a^{-q}N_-(\BH_{\a V})=
    \limsup\limits_{\a\to\w}\a^{-q}N_-(\BM_{\a G_V}),\\
    &\liminf\limits_{\a\to\w}\a^{-q}N_-(\BH_{\a V})=
    \liminf\limits_{\a\to\w}\a^{-q}N_-(\BM_{\a G_V}).
    \end{cases}
\end{equation*}
\end{prop}
This is an analog of statement (b) in Theorem 5.1 of the paper \cite{BL}. Its proof is basically the same, and we do not reproduce it here. In the same paper one finds also examples that illustrate the situation described by Proposition \ref{add}.
\section{Auxiliary material}\label{auxres}
The proof of \thmref{main} mainly follows the line worked out in \cite{BL} and \cite{S-dim2}. The same approach was used in \cite{LS}, and the material below, in part,
duplicates the contents of its Section 2.
We systematically use the variational description of the spectrum. In particular, we often define a self-adjoint operator via its corresponding Rayleigh quotient.
\subsection{Classes $\Sg_1,\Sg_1^\circ$ of compact operators}\label{clas}
If $\BT$ is a linear compact operator in a Hilbert space,
then, as usual, $\{s_j(\BT)\}$ stands for the sequence of its singular numbers, i.e., for the eigenvalues of the non-negative, self-adjoint operator $(\BT^*\BT)^{1/2}$. By $n_+(\vare,\BT)$ we denote the distribution function of the singular numbers,
\[n_+(\vare,\BT)=\#\{j:s_j>\vare\},\qquad \vare>0.\]

We say that $\BT$ belongs to the class $\Sg_1$ if and only if $\{s_j(\BT)\}\in\ell_{1,\w}$, and to the class $\Sg_1^\circ$ if and only if $\{s_j(\BT)\}\in\ell_{1,\w}^\circ$. These are linear, quasinormed
spaces with respect to the quasinorm $\|\BT\|_{1,\w}$ induced by this definition.  The space $\Sg_1$ is non-separable, and $\Sg_1^\circ$ is its
separable subspace in which the finite rank operators form a dense subset. Similarly to \eqref{Delta}, we define
the functionals
\begin{equation*}
    \D_1(\BT)=\D_1(\{s_j(\BT)\}),\qquad \d_1(\BT)=\d_1(\{s_j(\BT)\}).
\end{equation*}
Note that
\begin{equation*}
    \d_1(\BT)\le \D_1(\BT)\le\|\BT\|_{1,\w}.
\end{equation*}
See \cite{BSbook}, Section 11.6 for more detail about these spaces, and about the similar spaces $\Sg_q, \Sg_q^\circ$ for any $q>0$.

\subsection{Reduction of the main problem to compact operators}\label{comp}
 Let us introduce two subspaces in $C_0^\infty(\R^2)$:
\begin{equation*}
    \CF_0=\{f\in C_0^\infty: f(x)=\varf(r), \varf(1)=0\};\qquad \CF_1=\{f\in C_0^\infty: \int_0^{2\pi}f(r,\vart)d\vart=0,\ \forall r>0\}.
\end{equation*}
They are orthogonal to each other both in the $L_2$-metric and in the metric of the Dirichlet integral. The Hardy inequalities have a different form on $\CF_0$ and on $\CF_1$:
\begin{equation}\label{H0}
    \int\frac{|f(x)|^2}{|x|^2\ln^2|x|}dx\le\frac14\int|\nabla f(x)|^2dx,\qquad f\in \CF_0;
\end{equation}
\begin{equation}\label{H1}
    \int\frac{|f(x)|^2}{|x|^2}dx\le\int|\nabla f(x)|^2dx,\qquad f\in \CF_1.
\end{equation}

For proving \eqref{H0}, one substitutes $r=|x|=e^t$, and then applies the standard Hardy inequality in dimension 1.
The proof of \eqref{H1} is quite elementary, it can be found, e.g., in \cite{S-dim2}, or in \cite{BL}.

 Let us consider the completions $\CH^1_0, \CH^1_1$ of the spaces $\CF_0,\CF_1$ in the metric of the Dirichlet integral. It follows from the Hardy inequalities \eqref{H0}, \eqref{H1} that these are Hilbert function spaces, embedded into the weighted $L_2$, with the weights defined by these inequalities. Consider also their orthogonal sum
  \begin{equation}\label{ortsum}
    \CH^1=\CH_0^1\oplus\CH_1^1.
 \end{equation}
 An independent definition of this Hilbert space is
 \begin{equation*}
    \CH^1=\{f\in H^1_{\loc}(\R^2): \int_0^{2\pi}f(1,\vart)d\vart=0,\ |\nabla f|\in L_2(\R^2)\},
 \end{equation*}
 with the metric of the Dirichlet integral.

 We also define the spaces $H^1_0, H^1_1$ which are the completions of $\CF_0,\CF_1$ in $H^1(\R^2)$, and
 \[ \wt H^1=H_0^1\oplus H_1^1=\{f\in H^1(\R^2):\int_0^{2\pi}f(1,\vart)d\vart=0\}.\]
 This is a subspace in $H^1(\R^2)$ of codimension $1$.

  Finally, we need the spaces $\CG_0, \CG_1$ which are the completions of $\CF_0,\CF_1$ in the $L_2$-metric. Note that the condition $\varf(1)=0$, occuring in the description of $\CF_0$, disappears for general $f\in\CG_0$.

Suppose that $V\ge0$ is a measurable function, such that
 \begin{equation}\label{quaf}
    \bb_V[u]:=\int V|u|^2dx\le C\int|\nabla u|^2dx,\qquad\forall u\in\CH^1.
 \end{equation}
 Under the assumption \eqref{quaf} the quadratic form $\bb_V$ defines a bounded
 self-adjoint operator $\BB_V\ge0$ in $\CH^1$. If (and only if) this operator is compact, then, by the Birman-Schwinger principle, the quadratic form
 \begin{equation}\label{qf}
    \int(|\nabla u|^2-\a V|u|^2)dx
 \end{equation}
 with the form-domain $\wt H^1$
 is closed and bounded from below for each $\a>0$, the negative spectrum of the associated self-adjoint operator $\wt{\BH}_{\a V}$ on $L_2(\R^2)$ is finite, and the following equality for the number of its negative eigenvalues holds true:
 \begin{equation}\label{bshw}
    N_-(\wt{\BH}_{\a V})=n_+(\a^{-1},\BB_V),\qquad\forall\a>0.
 \end{equation}

 Now, let us withdraw the rank one condition $\int_0^{2\pi}u(1,\vart)d\vart=0$ from the description of the form-domain. Then the resulting quadratic form corresponds to the \sh operator $\BH_{\a V}$ as in \eqref{sh}. Hence,
 \[ N_-(\wt{\BH}_{\a V})\le N_-(\BH_{\a V})\le N_-(\wt{\BH}_{\a V})+1,\]
 and, by \eqref{bshw},
 \begin{equation*}
    n_+(\a^{-1},\BB_V)\le N_-(\BH_{\a V})\le n_+(\a^{-1},\BB_V)+1.
 \end{equation*}
 Thus, the study of the quantity $N_-(\BH_{\a V})$ for all $\a>0$ is reduced to the investigation of the "individual" operator $\BB_V$, which is
 nothing but the Birman-Schwinger operator for the family of operators in $L_2(\R^2)$
 associated with the family of quadratic forms in \eqref{qf}. Note that the Birman-Schwinger operator for the original family \eqref{sh} is ill-defined, since the completion of the space $H^1(\R^2)$ in the metric of the Dirichlet integral is not a space of functions on $\R^2$.
 \section{Proof of \thmref{main}}\label{ingr}
 \subsection{Decomposition of the quadratic form $\bb_V$.}
 Given a function $u\in\CH^1$, we agree to standardly denote its components in the decomposition \eqref{ortsum} by $\varf(r),v(r,\vart)$. Along with the quadratic form $\bb_V$, we consider its "parts" in the subspaces $\CH_0^1,\CH_1^1$:
 \begin{equation*}
    \bb_{V,0}[u]=\bb_V[\varf],\qquad \bb_{V,1}[u]=\bb_V[v].
 \end{equation*}
 Let $\BB_{V,j},\ j=0,1,$ stand for the corresponding self-adjoint operators in $\CH_j^1$. Using the orthogonal decomposition \eqref{ortsum}, we see that
 \begin{equation}\label{dec}
    \bb_V[u]=\bb_{V,0}[\varf]+\bb_{V,1}[v]+2\int V(x)\re(\varf(|x|)\overline{v(x)})dx.
\end{equation}
For the radial potentials the last term vanishes, and this considerably simplifies the reasoning,
see \cite{LS}. For the general potentials this is no more true. Still, the following inequality is always valid:
\begin{equation}\label{qform}
    \bb_V[u]\le 2(\bb_V[\varf]+\bb_V[v]),
\end{equation}
and it shows that for estimation of $\|\BB_V\|_{1,\w}$ it suffices to evaluate the
quasinorms in $\Sg_1$ of the operators $\BB_{V,0},\BB_{V,1}$ separately.

The estimation of $\|\BB_{V,0}\|_{1,\w}$ will be based upon the following result on the operators $\BF_G$ on real line, whose Rayleigh quotient is
\begin{equation}\label{ray1}
    \frac{\int_\R G(t)|\om(t)|^2dt}{\int_\R |\om'(t)|^2dt},\qquad \om(0)=0.
\end{equation}
Clearly, this is the Birman-Schwinger operator for the family $\BM_{\a G}$ given by \eqref{dim1}.

\begin{prop}\label{B0}
Let a function $G\in L_{1,loc}(\R),\ G\ge0,$ be given. Define the corresponding number sequence $\wh{\gz}(G)$ as in \eqref{sum}, and suppose that $\wh{\gz}(G)\in\ell_{1,\w}$.
Then the operator $\BF_G$ is well-defined, belongs to the class $\Sg_1$, and the estimate is satisfied,
\begin{equation}\label{est0}
    \|\BF_G\|_{1,\w}\le C\|\wh{\gz}(G)\|_{1,\w}.
\end{equation}
If $\wh{\gz}(G)\in\ell_{1,\w}^\circ$, then
$\wh{\gz}(G)\in\Sg_1^\circ$.
\end{prop}\vs

For the proof, see Section 4 in the paper \cite{BL}. There the operators on the half-line were considered, however the passage to the case of the whole line is straightforward, due to the
condition $\om(0)=0$ in \eqref{ray1}. In this respect, see also a discussion in \cite{LS}, Section 3.\vs

Now we turn to the operator $\BB_{V,1}$.
The estimation of its quasinorm in $\Sg_1$ uses a result that is a particular case (for $l=1$)
of Theorem 1.2 in the paper \cite{LN}. We present its equivalent formulation, more convenient for our purposes. Namely, we formulate it for the Birman-Schwinger operator, rather than for
the original \sh operator, as it was done in \cite{LN}.
\begin{prop}\label{lanet}
Let $V\ge0,\ V\in L_1(\R_+,\, L_p(\Sq))$, with some $p>1$. Then the operator $\wh{\BB}_V$, whose Rayleigh quotient is
\begin{equation}\label{ray2}
    \frac{\int V(x)|u|^2dx}{\int\left(|\nabla u|^2+|x|^{-2}|u|^2\right)dx},\qquad u\in \CH^1_1,
\end{equation}
belongs to the class $\Sg_1$, and
\begin{equation}\label{ray2-est}
    \|\wh{\BB}_V\|_{1,\w}\le C(p)\|V\|_{L_1(\R_+,\, L_p(\Sq))}.
\end{equation}
\end{prop}
We recall that the norm appearing in \eqref{ray2-est} was defined in \eqref{L1Lp}.

\subsection{Proof of \eqref{estim}}\label{main-est}
As it was explained in the previous subsection, we have to estimate the quasinorms of the operators $\BB_{V,0}, \BB_{V,1}$ in the space $\Sg_1$.

Consider first the operator $\BB_{V,0}$. The corresponding Rayleigh quotient is
\begin{equation}\label{missed}
    \frac{\int_{\R^2}V(r,\vart)|\varf(r)|^2rdrd\vart}
    {\int_{\R^2}|\varf'(r)|^2rdrd\vart}=
\frac{\int_0^\w V_{\rad}(r)|\varf(r)|^2rdr}{\int_0^\w|\varf'(r)|^2rdr}.
\end{equation}
The standard substitution $r=e^t,\ \varf(r)=\om(t);\ t\in\R$, reduces it to the form
\begin{equation*}
    \frac{\int_\R G_V(t)|\om(t)|^2dt}{\int_\R|\om'(t)|^2dt},\qquad \om(0)=0.
\end{equation*}
where the potential $G_V$ is given by \eqref{effpot}. Now, Proposition \ref{B0} applies, and we arrive at the estimate
\begin{equation*}
    \|\BB_{V,0}\|_{1,\w}\le C\|\wh{\gz}(G)\|_{1,\w}.
\end{equation*}

The Rayleigh quotient for the operator $\BB_{V,1}$ is given by
\begin{equation*}
    \frac{\int V(x)|u|^2dx}{\int|\nabla u|^2dx},\qquad u\in \CH^1_1.
\end{equation*}
Due to the Hardy inequality \eqref{H1}, on the subspace $\CH^1_1$ the norm of the Dirichlet integral is equivalent to the norm generated by the quadratic form in the denominator of \eqref{ray2}. Hence, the estimate \eqref{ray2-est} applies to this operator, with some other constant factor $C'(p)$. So, we have
\begin{equation}\label{ray3-est}
    \|\BB_{V,1}\|_{1,\w}\le C'(p)\|V\|_{L_1(\R_+,L_p(\Sq))}.
\end{equation}

The estimates \eqref{est0} and \eqref{ray3-est}, together with the inequality \eqref{qform}, imply the desired estimate \eqref{estim}.

\subsection{Proof of \eqref{as2}.}\label{prmain-as}
First of all, we are going to show that
\begin{equation}\label{as-B1}
    \lim_{\vare\to0}\left(\vare\ n_+(\vare,\BB_{V,1})\right)=\frac1{4\pi}\int_{\R^2}Vdx.
    \end{equation}

For $V\in C_0^\w(\R^2\setminus\{0\})$, Theorem 5.1 in \cite{BL} yields that $\N_-(\BH_{\a V})
\sim (4\pi)^{-1}\a\int Vdx$ as $\a\to\infty$. By the Birman-Schwinger principle, this is equivalent to $n_+(\vare,\BB_V)\sim (4\pi\vare)^{-1}\int Vdx$ as $\vare\to 0$. The spectrum of $\BB_{V,1}$ has the same asymptotic behavior, since for such potentials the subspace
$\CH^1_0$ does not contribute to the asymptotic coefficient.

Now, let $V\ge0$ be an arbitrary potential from $L_1(\R_+,L_p(\Sq))$. Then, approximating it by the
functions from $C_0^\w$ and taking into account the continuity of the asymptotic coefficients in the metric of $\Sg_1$ (see \cite{BSbook}, Theorem 11.6.6), we extend the formula to all such $V$. So, \eqref{as-B1} is established.\vs

Return to the study of the operator $\BB_V$.
Along with it, let us consider the direct orthogonal sum $\CB_V=\BB_{V,0}\oplus\BB_{V,1}$. Evidently,
\[ n_+(\vare,\CB_V)= n_+(\vare,\BB_{V,0})+n_+(\vare,\BB_{V,1}).\]
Hence, for justifying the asymptotic formulas \eqref{as2} it suffices to show that the off-diagonal term in \eqref{dec} generates an operator of the class $\Sg_1^\circ$. To this end, we first of all note that
\begin{equation}\label{reduc}
    \int V\re(\varf\overline{v})dx=\int V_{\nrad}\re(\varf\overline{v})dx,
\end{equation}
since $v$ is orthogonal (in $L_2$) to all functions depending only on $|x|$.

Suppose now that the function $V_{\nrad}$ has a compact support in $\R^2\setminus\{0\}$. Then the integral in the right-hand side of \eqref{reduc} is actually taken over some annulus $a\le r\le a^{-1},\ a<1$. Hence,
\begin{gather*}
2\bigl|\int V_{\nrad} \re(\varf\overline{v})dx\bigr|
\le\d\int_a^{a^{-1}}rdr\int_\Sq|V_{\nrad}(r,\vart)|
|v(r,\vart)|^2d\vart\\ +
\d^{-1}\int_a^{a^{-1}}rdr\int_{\Sq}|V_{\nrad}(r,\vart)|
|\varf(r)|^2d\vart.
\end{gather*}
The first term on the right generates an operator on $\CH^1_1$, say, $\BT_1$, to which the estimate \eqref{ray3-est} applies, and it gives
\[ \|\BT_1\|_{1,\w}\le C'(p)\d.\]
The second term generates an operator on $\CH^1_0$, say, $\BT_0$. Its Rayleigh quotient is of the type
\eqref{missed} but with the integration over a compact subset in $(0,\infty)$. It follows that the spectrum of $\BT_0$ obeys Weyl's asymptotic law, $\l_j(\BT_0)\asymp c\j^{-2}$, and hence, $\BT_0\in\Sg_1^\circ$. Taking $\d$ arbitrarily small, we conclude that the asymptotics \eqref{as2} is satisfied in the case  where $V_{\nrad}$ is compactly supported.

Finally, we approximate the function $V_{\nrad}$ by compactly supported functions in the metric \eqref{L1Lp}, and again apply Theorem 11.6.5 from the book \cite{BSbook}. This extends the asymptotic formula \eqref{as2} to all potentials, that meet the conditions of \thmref{main}, and thus, concludes the proof.



\begin{thebibliography}{22}

\bibitem{BL} Birman, M.Sh., Laptev, A.: {\it The negative discrete spectrum
of a two-dimensional Schr\"odinger operator,} Comm. Pure Appl. Math. {\bf 49},  no. 9, 967--997  (1996).

\bibitem{BSbook} Birman, M.Sh., Solomyak, M.: {\it Spectral theory of
selfadjoint operators in Hilbert space}, D. Reidel Publishing Co.,
Dordrecht, 1987.

\bibitem{cha} Chadan, K., Khuri, N. N., Martin, A., Wu, Tai Tsun: {\it Bound states
in one and two spatial dimensions,}  J. Math. Phys. {\bf 44},  no. 2, 406--422 (2003).

\bibitem{GN} Grigor'yan, A., Nadirashvili, N.:
{\it Negative eigenvalues of two-dimensional Schr\"odinger operators,} arXiv:1112.4986

\bibitem{Lap} Laptev, A.: {\it The negative spectrum of the class of
two-dimensional Schr\"odinger operators with
potentials that depend on the radius,} (Russian)
Funktsional. Anal. i Prilozhen. {\bf 34}, no.
4, 85--87  (2000);  translation in  Funct. Anal. Appl.
{\bf 34}, no. 4, 305--307  (2000).

\bibitem{LN} Laptev, A., Netrusov, Yu.:
{\it On the negative eigenvalues of a class of Schr\"odinger operators},  Differential operators and spectral theory, 173--186,
Amer. Math. Soc. Transl. Ser. 2, {\bf 189}, Amer. Math. Soc., Providence, RI, 1999.

\bibitem{LS} Laptev, A., Solomyak, M.: {\it On the negative spectrum of the two-dimensional Schr\"odinger operator with radial potential}, Comm. Math. Phys., in press.

\bibitem{mv} Molchanov, S., Vainberg, B.: {\it On negative eigenvalues of low-dimensional Schr\"odinger operators}, arXiv:1105.0937.

\bibitem{S-dim2} Solomyak, M.: {\it Piecewise-polynomial approximation of functions from
$H^l((0,1)^d),\ 2l=d$, and applications to the
spectral theory of the Schr\"odinger operator}, Israel
J. Math.  {\bf 86}, no. 1-3, 253--275  (1994).

\end{thebibliography}
\end{document}